\documentclass[14pt]{amsart}

\usepackage{amsmath, amssymb, amsthm}

\usepackage[figuresright]{rotating}

\usepackage{fullpage}

\usepackage{color}
\usepackage{etex}
\usepackage{dsfont}
\usepackage[matrix, arrow, curve, frame, arc]{xy}
\usepackage{pgf}
\usepackage{tikz}
\usepackage{amscd}
\usetikzlibrary{arrows,shapes,automata,backgrounds,petri,calc}

\newcommand{\tikzAngleOfLine}{\tikz@AngleOfLine}
\def\tikz@AngleOfLine(#1)(#2)#3{%
\pgfmathanglebetweenpoints{%
\pgfpointanchor{#1}{center}}{%
\pgfpointanchor{#2}{center}}
\pgfmathsetmacro{#3}{\pgfmathresult}%
}

\newcommand{\bN}{\mathbb{N}}

\newcommand{\wt}{\widetilde}
\newcommand{\ba}{\bar{\alpha}}
\newcommand{\La}{\Lambda}

\newcommand{\ve}{\varepsilon}

\newtheorem{theorem}{Theorem}[section] 
\newtheorem{lemma}[theorem]{Lemma}   
\newtheorem{corollary}[theorem]{Corollary}
\newtheorem{proposition}[theorem]{Proposition}

\newtheorem{main-theorem}[theorem]{Theorem}

\newtheorem*{problem*}{Problem}

\theoremstyle{definition}

\newtheorem*{question*}{Question}

\begin{document}

\baselineskip=14pt
\title[Spherical algebras]{A note on spherical algebras$^*$\footnote{\tiny $^*$ The research has been supported 
from the research grant no. 2023/51/D/ST1/01214 of the Polish National Science Center}} 

\author[K. Erdmann]{Karin Erdmann}
\address[Karin Erdmann]{Mathematical Institute, University of Oxford, UK}
\email{erdmann@maths.ox.ac.uk}

\author[A. Skowyrski]{Adam Skowyrski}
\address[Adam Skowyrski]{Faculty of Mathematics and Computer Science, 
Nicolaus Copernicus University, Chopina 12/18, 87-100 Toru\'n, Poland}
\email{skowyr@mat.umk.pl}


\begin{abstract} We classify tame symmetric algebras of period four which are
	closely related to the spherical algebras introduced in \cite{WSA-GV}. This note provides a classification 
	in the special case which naturally appears, when dealing with biregular Gabriel quivers. 
\end{abstract}

\maketitle

\section{\bf Introduction}

This is part of the project of classifying the tame symmetric algebras for which all simple modules have period four, 
called TSP4 algebras. TSP4 algebras include most weighted surface algebras (WSA's), and various other algebras, 
introduced and studied  in \cite{WSA}, \cite{WSA-GV}, \cite{HTA}, \cite{HSA}, \cite{HSS}, \cite{SS}. \medskip  

The algebras are tame, which implies that the Gabriel quiver (that is, the ordinary quiver) of the algebra cannot 
have too many arrows starting or ending at any vertex. As well the quiver does not have sources or sinks since the
algebra is symmetric. TSP4 algebras with $2$-regular Gabriel quivers (that is, two arrows start and end at each vertex) 
were classified in \cite{GQT}. The present stage is to obtain a classification when the Gabriel quiver is biregular, i.e. 
the numbers of arrows starting and ending at given vertex are equal (regular), and do not exceed $2$. Equivalently, 
each vertex is either a $2$-vertex (two arrows start and two arrows end) or is a $1$-vertex (one arrow starts and 
one arrow ends).

We note that the classification in the $2$-regular case \cite{GQT} was based on the first (simplified) definition 
of WSA's \cite{WSA}, where the weights were restricted to get a triangulated structure on the Gabriel quiver of 
the algebra. In this case, for any TSP4 algebra $\La=KQ/I$, its Gabriel quiver $Q$ is a glueing of a finite number 
of three types of blocks: 
$$\mbox{I: }
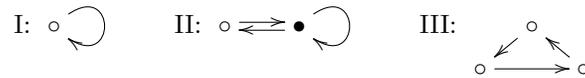
\xymatrix@C=0.3cm@R=0.2cm{\ar@(ru, rd)[]\circ} \qquad\qquad \mbox{II: }
\xymatrix@C=0.6cm@R=0.2cm{\circ\ar@<.35ex>[r]&\ar@<.35ex>[l]\bullet\ar@(ru,rd)[] } \qquad\qquad  \mbox{III: } 
\xymatrix@C=0.3cm@R=0.2cm{&\circ\ar[ld]_{}&\\ \circ\ar[rr]_{}&&\circ\ar[lu]_{}}$$ 
where by a glueing we mean that each vertex $\circ$ is glued with exactly one vertex $\circ$ in a different block. 
In general \cite{WSA-GV}, a weighted surface algebra $\La$ is given as a quotient $\La=KQ/I$, where $Q$ is a glueing 
of a finite number of blocks of types I-III, or equivalently, $Q$ is a triangulation quiver (see also Section 
\ref{2.1}). The ideal $I$ depends on the chosen weights, and it can happen that some of the arrows in $Q$ are not 
arrows of the Gabriel quiver $Q_\La$ of $\La$ (called virtual arrows). Following \cite{WSA-GV}, we conclude that 
the Gabriel quiver of any weighted surface algebra $\La$ is a glueing of finite number of block of types I-III and 
two new types of blocks of the form 
$$
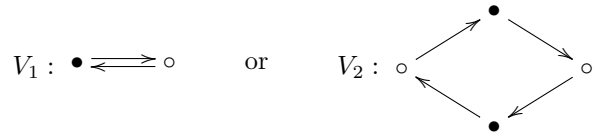
\xymatrix@R=0.4cm{\\ V_1: \\ }\xymatrix@R=0.4cm{&\\ \bullet\ar@<.35ex>[r]&\ar@<.35ex>[l]\circ}\qquad 
\xymatrix@R=0.4cm{\\ \mbox{or}\\} \qquad 
\xymatrix@R=0.4cm{\\ V_2: \\ }\xymatrix@R=0.4cm{&\bullet \ar[rd]& \\ \circ \ar[ru] && \circ \ar[ld]\\ & \bullet\ar[lu]&}$$ 
The blocks of types $V_1,V_2$ are obtained, respectively, from a block of type II by removing a virtual loop, 
or from a glueing of two triangles by removing a $2$-cycle of virtual arrows (see also the introduction of 
\cite{EHS2}). \medskip 

The first step towards the classification in the biregular case was to understand the position of $1$-vertices 
in the Gabriel quiver of the algebra. In \cite{EHS2} we managed to show that for any TSP4 algebra $\La=KQ/I$ with  
biregular Gabriel quiver $Q$, every $1$-vertex of $Q$ is contained in a block of type $V_1$ or $V_2$. This 
confirmes the conjecture that all such algebras should have the same Gabriel quivers as WSA's, but this 
requires a bit more work. In a current project, we try to prove that every TSP4 algebra with biregular quiver 
is a WSA (with some exceptions as in \cite{GQT}), or is the so called higher spherical algebra HSA \cite{HSA}, 
which is not a WSA  but it has the spherical quiver (a glueing of two blocks of type V$_2$, discussed 
below). It turned out that the proof naturally splits into two cases, where the first case embraces all 
algebras with Gabriel quivers different from the spherical quiver (or the so called triangle quiver, which is 
not considered in our paper; see \cite[Section 6.5]{EHS3} for details in this case). 

In this note, we cover the second case, that is, we classify TSP4 algebras having  the spherical quiver. 
More precisely, we prove that such algebras are either WSA's or HSA's. A parallel paper \cite{EHS3} is devoted 
to a classification in the first case, so this article should be seen as a missing case to complete the 
classification in general. This is the one of two reasons for dealing with this quivers separately. \medskip 
 
Here, we distinguish the following two types of `spherical' quivers, denoted by $Q^S$ and $Q^{S'}$, respectively.  
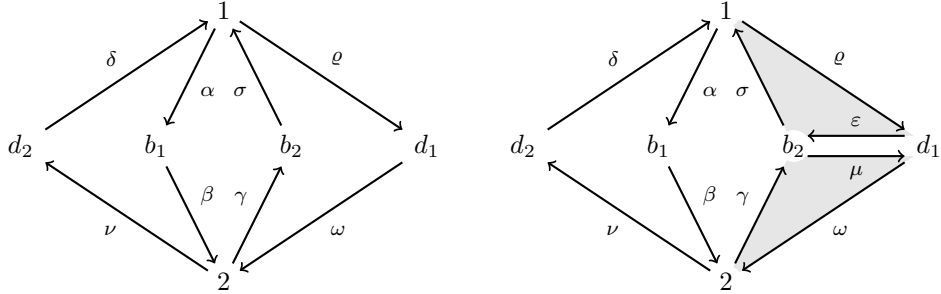
\[ \begin{tikzpicture}
[->,scale=.9]
\coordinate (1) at (0,2);
\coordinate (2) at (-1,0);
\coordinate (2u) at (-.925,.15);
\coordinate (2d) at (-.925,-.15);
\coordinate (3) at (0,-2);
\coordinate (4) at (1,0);
\coordinate (4u) at (.925,.15);
\coordinate (4d) at (.925,-.15);
\coordinate (5) at (-3,0);
\coordinate (5u) at (-2.775,.15);
\coordinate (5d) at (-2.775,-.15);
\coordinate (6) at (3,0);
\coordinate (6u) at (2.775,.15);
\coordinate (6d) at (2.775,-.15);
\node [fill=white,circle,minimum size=4.5] (1) at (0,2) {\ \quad};
\node [fill=white,circle,minimum size=4.5] (2) at (-1,0) {\ \quad};
\node [fill=white,circle,minimum size=4.5] (3) at (0,-2) {\ \quad};
\node [fill=white,circle,minimum size=4.5] (4) at (1,0) {\ \quad};
\node [fill=white,circle,minimum size=4.5] (5) at (-3,0) {\ \quad};
\node [fill=white,circle,minimum size=4.5] (6) at (3,0) {\ \quad};
\node (1) at (0,2) {1};
\node (2) at (-1,0) {$b_1$};
\node (2u) at (-1,0.15) {\ \quad};
\node (2d) at (-1,-0.15) {\ \quad};
\node (3) at (0,-2) {2};
\node (4) at (1,0) {$b_2$};
\node (4u) at (1,0.15) {\ \quad};
\node (4d) at (1,-0.15) {\ \quad};
\node (5) at (-3,0) {$d_2$};
\node (5u) at (-2.775,0.15) {};
\node (5d) at (-2.775,-0.15) {};
\node (6) at (3,0) {$d_1$};
\node (6u) at (2.775,0.15) {};
\node (6d) at (2.775,-0.15) {};
\draw[thick,->]
(1) edge node[below right]{\footnotesize$\alpha$} (2)
(5u) edge node[above left]{\footnotesize$\delta$} (1)
(2) edge node[above right]{\footnotesize$\beta$} (3)
(3) edge node[below left]{\footnotesize$\nu$} (5d)
(1) edge node[above right]{\footnotesize$\varrho$} (6u)
(4) edge node[below left]{\footnotesize$\sigma$} (1)
(6d) edge node[below right]{\footnotesize$\omega$} (3)
(3) edge node[above left]{\footnotesize$\gamma$} (4)
;
\end{tikzpicture}\qquad
\begin{tikzpicture}
[->,scale=.9]
\coordinate (1) at (0,2);
\coordinate (2) at (-1,0);
\coordinate (2u) at (-.925,.15);
\coordinate (2d) at (-.925,-.15);
\coordinate (3) at (0,-2);
\coordinate (4) at (1,0);
\coordinate (4u) at (.925,.15);
\coordinate (4d) at (.925,-.15);
\coordinate (5) at (-3,0);
\coordinate (5u) at (-2.775,.15);
\coordinate (5d) at (-2.775,-.15);
\coordinate (6) at (3,0);
\coordinate (6u) at (2.775,.15);
\coordinate (6d) at (2.775,-.15);
\fill[fill=gray!20] (1) -- (4u) -- (6u) -- cycle;
\fill[fill=gray!20] (3) -- (6d) -- (4d) -- cycle;
\node [fill=white,circle,minimum size=4.5] (1) at (0,2) {\ \quad};
\node [fill=white,circle,minimum size=4.5] (2) at (-1,0) {\ \quad};
\node [fill=white,circle,minimum size=4.5] (3) at (0,-2) {\ \quad};
\node [fill=white,circle,minimum size=4.5] (4) at (1,0) {\ \quad};
\node [fill=white,circle,minimum size=4.5] (5) at (-3,0) {\ \quad};
\node [fill=white,circle,minimum size=4.5] (6) at (3,0) {\ \quad};
\node (1) at (0,2) {1};
\node (2) at (-1,0) {$b_1$};
\node (2u) at (-1,0.15) {\ \quad};
\node (2d) at (-1,-0.15) {\ \quad};
\node (3) at (0,-2) {2};
\node (4) at (1,0) {$b_2$};
\node (4u) at (1,0.15) {\ \quad};
\node (4d) at (1,-0.15) {\ \quad};
\node (5) at (-3,0) {$d_2$};
\node (5u) at (-2.775,0.15) {};
\node (5d) at (-2.775,-0.15) {};
\node (6) at (3,0) {$d_1$};
\node (6u) at (2.775,0.15) {};
\node (6d) at (2.775,-0.15) {};
\draw[thick,->]
(1) edge node[below right]{\footnotesize$\alpha$} (2)
(5u) edge node[above left]{\footnotesize$\delta$} (1)
(2) edge node[above right]{\footnotesize$\beta$} (3)
(3) edge node[below left]{\footnotesize$\nu$} (5d)
(1) edge node[above right]{\footnotesize$\varrho$} (6u)
(6u) edge node[above]{\footnotesize$\varepsilon$} (4u)
(4) edge node[below left]{\footnotesize$\sigma$} (1)
(4d) edge node[below]{\footnotesize$\mu$} (6d)
(6d) edge node[below right]{\footnotesize$\omega$} (3)
(3) edge node[above left]{\footnotesize$\gamma$} (4)
;
\end{tikzpicture}
\]
The first one is a glueing of two blocks of type V$_2$ and the second is a glueing of one block of type 
V$_2$ and two triangles (blocks of type III). The quiver $Q^S$ will be called the spherical quiver, whereas 
$Q^{S'}$, the almost spherical. Both are subquivers of the following quiver $Q$ defining the spherical algebras as in 
\cite[see Example 3.6]{WSA-GV}. 

\[
\begin{tikzpicture}
[->,scale=.9]
\coordinate (1) at (0,2);
\coordinate (2) at (-1,0);
\coordinate (2u) at (-.925,.15);
\coordinate (2d) at (-.925,-.15);
\coordinate (3) at (0,-2);
\coordinate (4) at (1,0);
\coordinate (4u) at (.925,.15);
\coordinate (4d) at (.925,-.15);
\coordinate (5) at (-3,0);
\coordinate (5u) at (-2.775,.15);
\coordinate (5d) at (-2.775,-.15);
\coordinate (6) at (3,0);
\coordinate (6u) at (2.775,.15);
\coordinate (6d) at (2.775,-.15);
\fill[fill=gray!20] (1) -- (5u) -- (2u) -- cycle;
\fill[fill=gray!20] (1) -- (4u) -- (6u) -- cycle;
\fill[fill=gray!20] (2d) -- (5d) -- (3) -- cycle;
\fill[fill=gray!20] (3) -- (6d) -- (4d) -- cycle;
\node [fill=white,circle,minimum size=4.5] (1) at (0,2) {\ \quad};
\node [fill=white,circle,minimum size=4.5] (2) at (-1,0) {\ \quad};
\node [fill=white,circle,minimum size=4.5] (3) at (0,-2) {\ \quad};
\node [fill=white,circle,minimum size=4.5] (4) at (1,0) {\ \quad};
\node [fill=white,circle,minimum size=4.5] (5) at (-3,0) {\ \quad};
\node [fill=white,circle,minimum size=4.5] (6) at (3,0) {\ \quad};
\node (1) at (0,2) {1};
\node (2) at (-1,0) {$b_1$};
\node (2u) at (-1,0.15) {\ \quad};
\node (2d) at (-1,-0.15) {\ \quad};
\node (3) at (0,-2) {2};
\node (4) at (1,0) {$b_2$};
\node (4u) at (1,0.15) {\ \quad};
\node (4d) at (1,-0.15) {\ \quad};
\node (5) at (-3,0) {$d_2$};
\node (5u) at (-2.775,0.15) {};
\node (5d) at (-2.775,-0.15) {};
\node (6) at (3,0) {$d_1$};
\node (6u) at (2.775,0.15) {};
\node (6d) at (2.775,-0.15) {};
\draw[thick,->]
(1) edge node[below right]{\footnotesize$\alpha$} (2)
(2u) edge node[above]{\footnotesize$\xi$} (5u)
(5u) edge node[above left]{\footnotesize$\delta$} (1)
(5d) edge node[below]{\footnotesize$\eta$} (2d)
(2) edge node[above right]{\footnotesize$\beta$} (3)
(3) edge node[below left]{\footnotesize$\nu$} (5d)
(1) edge node[above right]{\footnotesize$\varrho$} (6u)
(6u) edge node[above]{\footnotesize$\varepsilon$} (4u)
(4) edge node[below left]{\footnotesize$\sigma$} (1)
(4d) edge node[below]{\footnotesize$\mu$} (6d)
(6d) edge node[below right]{\footnotesize$\omega$} (3)
(3) edge node[above left]{\footnotesize$\gamma$} (4)
;
\end{tikzpicture}
\]

We recall that $Q$ is a triangulation quiver, with permutation $f$ of the arrows given by four $3$-cycles 
(rotating arrows in the four shaded triangles). Then the associated permutation $g=\bar{f}$, where $\overline{(-)}$ 
is the involution of the set of arrows (for details, see Section \ref{2.1}), has two orbits of length 
four: 
$$(\alpha \ \beta \ \gamma \ \sigma)\mbox{ and }(\rho \ \omega \ \nu \ \delta)$$ 
and two orbits of length two: $(\xi \ \eta)$ and $(\ve \ \mu)$. Let $m_\bullet$ and $c_\bullet$ be 
the weight and parameter functions for the triangulation quiver $(Q,f)$, and $\La=KQ/I$ be the 
associated weighted surface algebra $\La=\La(Q,f,m_\bullet,c_\bullet)$. Then, depending on the weights 
$n=m_\xi$ and $n'=m_{\xi}$, the Gabriel quiver of $\La$ and the relations may differ. For a precise 
description of relations in $I$ in general case, we refer to Section \ref{2.1}. \smallskip 

We will not consider WSA's with $n,n'\geqslant 2$, since then $Q_\La=Q$ is $2$-regular, and all TSP4 
algebras with this Gabriel quiver were classified in \cite{GQT}. In the case when $n=n'=1$, arrows 
$\xi,\eta,\ve,\mu$ are all virtual, and the Gabriel quiver of $\La$ is $Q^S$. Then $\La$ is 
given by $Q^S$, with the admissible ideal generated by the following commutativity relations 
\begin{align*}
(S1a) &\ \beta\nu\delta=a(\beta\gamma\sigma\alpha)^{m-1}\beta\gamma\sigma, &&& 
(S1b) &\ \nu\delta\alpha=a(\gamma\sigma\alpha\beta)^{m-1}\gamma\sigma\alpha, & \\ 
(S2a) &\ \delta\alpha\beta=b(\delta\rho\omega\nu)^{m'-1}\delta\rho\omega, &&& 
(S2b) &\ \alpha\beta\nu=b(\rho\omega\nu\delta)^{m'-1}\rho\omega\nu, & \\
(S3a) &\ \sigma\rho\omega=a(\sigma\alpha\beta\gamma)^{m-1}\sigma\alpha\beta,&&& 
(S3b) &\ \rho\omega\gamma=a(\alpha\beta\gamma\sigma)^{m-1}\alpha\beta\gamma, & \\ 
(S4a) &\ \omega\gamma\sigma=b(\omega\nu\delta\rho)^{m'-1}\omega\nu\delta,&&& 
(S4b) &\ \gamma\sigma\rho=b(\nu\delta\rho\omega)^{m'-1}\nu\delta\rho, & 
\end{align*} 
and the following zero relations  
\begin{align*} 
(Z1) \ \alpha\beta\nu\delta\alpha=0, && (Z2)\ \nu\delta\alpha\beta\nu=0, && (Z3)\ \beta\nu\delta\rho=0,   
&& (Z4) \ \delta\alpha\beta\gamma=0,\quad & (Z5) \ \beta\nu\delta\alpha\beta=0,  & \\ 
(Z6) \ \delta\alpha\beta\nu\delta=0, && (Z7) \ \beta\gamma\sigma\rho=0,&& (Z8) \ 
\delta\rho\omega\gamma=0,&& (Z9) \ \gamma\sigma\rho\omega\gamma=0, \quad & (Z10) \ \rho\omega\gamma\sigma\rho=0, & \\ 
(Z11) \ \omega\gamma\sigma\alpha=0, && (Z12) \ \sigma\rho\omega\nu=0,&& 
(Z13) \ \sigma\rho\omega\gamma\sigma=0, && (Z14) \ \omega\gamma\sigma\rho\omega=0, \quad & 
(Z15) \ \sigma\alpha\beta\nu=0, & \\ 
&&&&(Z16) \ \omega\nu\delta\alpha=0, &&&& 
\end{align*} 
where $m=m_\alpha$, $m'=m_\rho$, $a=c_{\alpha}$ and $b=c_\rho$  (we can assume that other 
parameters are $1$). Such algebras will be sometimes called spherical algebras.   

If one of $n,n'$ is $1$, say $n=1$ and $n'\geqslant 2$, then the Gabriel quiver of $\La$ is the 
almost spherical quiver, with two virtual arrows $\xi,\eta$, and the admissible ideal for 
$\La$ is again generated by the commutativity relations and zero relations. We have the following 
commutativity relations: 
\begin{align*}
(S'1a) &\ \beta\nu\delta=a(\beta\gamma\sigma\alpha)^{m-1}\beta\gamma\sigma, &&& 
(S'1b) &\ \nu\delta\alpha=a(\gamma\sigma\alpha\beta)^{m-1}\gamma\sigma\alpha, & \\ 
(S'2a) &\ \delta\alpha\beta=b(\delta\rho\omega\nu)^{m'-1}\delta\rho\omega, &&& 
(S'2b) &\ \alpha\beta\nu=b(\rho\omega\nu\delta)^{m'-1}\rho\omega\nu, & \\
(S'3a) &\ \mu\omega=a(\sigma\alpha\beta\gamma)^{m-1}\sigma\alpha\beta,&&& 
(S'3b) &\ \rho\ve=a(\alpha\beta\gamma\sigma)^{m-1}\alpha\beta\gamma, & \\ 
(S'4a) &\ \ve\sigma=b(\omega\nu\delta\rho)^{m'-1}\omega\nu\delta,&&& 
(S'4b) &\ \gamma\mu=b(\nu\delta\rho\omega)^{m'-1}\nu\delta\rho. & 
\end{align*} 
where $m,m',a,b$ are defined as for spherical algebras. The zero relations in this case are $(Z'1)-(Z'20)$, where 
$(Z'1)-(Z'6)$ coincide with $(Z1)$-$(Z6)$, $(Z'7)$-$(Z'16)$ are obtained from the corresponding relations 
$(Z7)$-$(Z16)$ by replacing paths $\sigma\rho$ and $\omega\gamma$, by arrows $\mu$ and $\ve$, respectively, 
and there are additional relations $(Z'17)$-$(Z'20)$: $\sigma\rho\omega=0,\ \omega\gamma\sigma=0,\ 
\gamma\sigma\rho=0, \ \rho\omega\gamma =0$

The other important property, which makes this study interesting, is that on the spherical quiver $Q^S$ 
there are other structures of TSP4 algebras (besides the spherical algebras). As we mentioned earlier, 
there exist so called Higher Spherical Algebras (HSA), which are TSP4 algebras, but not WSA's. These are 
given by the spherical quiver $Q^S$ (which is biregular, but not 2-regular), and a modified set of relations 
(see Section \ref{5} for details). In \cite{HSA} the HSA was constructed via a tilting complex, showing it is
derived equivalent to algebras we call Higher Tetrahedral Algebras (HTA's), but it did not 
consider classifications of all TSP4 algebras with this quiver (see \cite[Section 5]{HSS}, for another 
approach to HSA's using mutations). Furthermore, this will illustrate the methods which are used in the 
general case. \medskip 

Our first result is the following. 

\medskip
\begin{theorem}\label{thm:1.1} Assume $\La$ is a TSP4 algebra with Gabriel quiver $Q^S$. Then one of the following holds:\\
(a) \  $\La$ is isomorphic to a weighted surface algebra $KQ/I$ (different from the singular spherical algebra), 
where $Q$ is the spherical quiver with four virtual arrrow $\xi, \eta, \ve, \mu$.\\
(b) \ $\La$ is isomorphic to a Higher Spherical Algebra $S(m, \lambda)$ with $m>1$. 
\end{theorem}

We complete the picture by proving the following theorem.  

\begin{theorem}\label{thm:1.2} Assume $\La$ is a TSP4 algebra with Gabriel quiver $Q^{S'}$. Then $\La$ is 
isomorphic to a weighted surface algebra $KQ/I$ where $Q$ is the spherical quiver with two virtual arrows $\xi,  \eta$. 
These algebras belong to the same family of weighted surface algebras as those in part (a) of the previous Theorem. \end{theorem} \bigskip

The quiver $Q^S$  is glued from two blocks of type V$_2$, up to labelling
we may take these blocks as
$$(1, \ b_1, \ 2, \ d_2) \ \ \mbox{and} \ (1, \ d_1, \ 2, \ b_2). 
$$
The quiver $Q^{S'}$ is glued from one block of type V$_2$, that is, $(1, \ b_1, \ 2, \ d_2)$, and two 
blocks of type  III (the shaded triangles).

\bigskip

Let us briefly describe the strategy for the classification. We start with a TSP4 algebra $\La$ having $Q^S$ or $Q^{S'}$ as 
the Gabriel quiver, and first determine minimal relations satisfied by $\La$. Then we describe bases for the 
indecomposable projective $\La$-modules, and hence their dimensions. In each case, there will be a unique WSA, 
or HSA respectively, $\hat{\La}$, satisfying the same relations. We will then define a surjective algebra 
homomorphism $\psi: \hat{\La} \to \La$ and show that it is an isomorphism.

The paper is organized as follows. In Section \ref{sec:2} we describe WSA's, and discuss minimal relations 
involving paths of length $\leq 3$ coming from the fact that the algebra is tame. Then we describe exact sequences 
determining the syzygies of the simple modules, and explain how to use these to deal explicitly with minimal relations 
of the algebra.

Section \ref{3} determines explicit minimal relations and spanning sets of the projective modules at $2$-vertices. 
There are two scalar parameters occuring in the minimal relations. Section \ref{4} deals with cases when at least one of these 
scalars is zero, and proves Theorem \ref{thm:1.1} in this case. In Section \ref{5} we determine algebras when both 
scalars are non-zero, and prove that they are HSA's denoted by $S(\lambda, m)$ for $m\geq 2$, or WSA's when 
$m=1$. The last section finishes the proof of Theorem \ref{thm:1.2}.

{\bf Acknowledgements.} Both authors thank the MFO Oberwolfach program ``Research in Pairs", the Mathematical 
Institute in Oxford and the Faculty of Mathematics and Computer Science of the Nicolaus Copernicus University 
for their hospitality. The second named author has been supported from the research grant no. 
2023/51/D/ST1/01214 of the Polish National Science Center.

\bigskip

\section{\bf Preliminaries}\label{sec:2}

For general background we refer to \cite{ASS}.
We assume throughout that algebras are basic and symmetric. 
Recall that the Gabriel quiver $Q_{\La}$ (that is, the ordinary quiver,) of $\La$
has vertices $Q_0$ in bijection with the simple modules, and 
arrows $i\to j$ correspond to a basis of $e_i(J/J^2)e_j$ for $i, j\in Q_0$.

Assume $\La$ is basic and symmetric, with Gabriel quiver $Q_{\La}$. Then $\La$ is a factor algebra of $KQ_{\La}$. A path in $Q_{\La}$
starting with an arrow $\tau: i\to j$, which is non-zero in $\La$, and is maximal
with this property, must  end in $i$ (since $\La$ is symmetric). We may choose
such a path, $W_{\tau}$, its image in $\La$ spans the socle of
$e_i\La$. We will  take such a
path which is  the  power of a
        cyclic path $X_{\tau}$ from $i$ to $i$ which is not the power of a
        shorter path. 
The algebra is symmetric, and therefore all rotations of such a path 
$W_{\tau}$ are also non-zero in $\La$ and belong to its socle.
We will use the following  convention:
        If $p$ and $q$ are paths of the same length $r$ satisfying  $p - cq \in J^{r+1}$ for some $0\neq c\in K$, we
write $p\equiv q$. 

For a path $\rho$ in $Q$, we will use notation $\rho\prec I$, if the path $\rho$ is a summand of some relation 
in a minimal set of generators for $I$.

\subsection{Weighted surface algebras}\label{2.1}

We recall the definition from \cite{WSA-GV-corr}. A weighted surface algebra is defined in terms of a presentation $\La = KQ/I = \La(Q, f, m_{\bullet}, c_{\bullet})$, the data 
$(Q, f, m_{\bullet}, c_{\bullet})$ defining the algebra 
are a triangulation quiver $(Q, f)$,  a weight function $m_{\bullet}$ and
a parameter  function  $c_{\bullet}$.

The quiver $Q$  is $2$-regular if at each vertex
two arrows start and two arrows end. Such a quiver has an involution $\alpha\to \ba$ on the arrows, where $\ba \neq \alpha$ are the two arrows starting at the same vertex.
A triangulation quiver is a pair $(Q, f)$ where $Q$ is a (finite) 2-regular quiver, and $f$ is a permutation of the arrows such that $f(\alpha)$ starts at the
vertex at which $\alpha$ ends, and moreover $f^3 = {\rm id}$. The permutation $f$ uniquely determines a permutation $g$ such that $g(\alpha) := \overline{f(\alpha)}$.
The  weight function
is a function $m_{\bullet}: Q_1\to \bN^*$ constant on $g$-cycles, and a parameter function is a function $c_{\bullet}: Q_1\to K^*$, also constant on $g$-cycles.
Let $n_{\alpha}$ be the length of the $g$-orbit of $\alpha$. We assume that $m_{\alpha}n_{\alpha} \geq 2$ for all arrows $\alpha$.
If $m_{\alpha}n_{\alpha}=2$ then $\alpha$ is a {\it virtual}  arrow.

For an arrow $\alpha$, we denote by $A_{\alpha}$ the path of length $m_{\alpha}n_{\alpha}-1$ along the $g$-cycle of $\alpha$.
In general, one assumes $m_{\alpha}n_{\alpha}\geq 3$ if $\ba$ is virtual and not a loop, and $m_{\alpha}n_{\alpha}\geq 4$ if
$\ba$ is a virtual loop.
With these data,
 $\La = KQ/I$ is a weighted surface algebra
if $I$ is generated by the following relations:

\begin{itemize}\item [(i)] $\alpha f\alpha- c_{\ba}A_{\ba}$ for all arrows $\alpha$ of $Q$,
                \item [(ii)] $\alpha f(\alpha)g(f(\alpha))$ for all  arrows $\alpha$ unless $f^2(\alpha)$ is virtual, or unless $f(\ba)$ is virtual and
$m_{\ba}=1, n_{\ba} = 3$.
                \item[(iii)] $\alpha g(\alpha)f(g(\alpha))$ for all
                        arrows $\alpha$ of $Q$ unless $f(\alpha)$ is virtual, or
                        unless $f^2(\alpha)$ is virtual, and $m_{f(\alpha)}=1$, $n_{f(\alpha)} = 3$.
        \end{itemize}

Note that $m_{\ba}n_{\ba}= 2$ (that is, $\ba$ is virtual,)  is equivalent with $A_{\ba} = \ba$. The relations show that then 
$\ba$ lies in $J^2$, and so is not part of the Gabriel quiver of $\La$. In fact, $Q$ is the disjoint union
of $Q_{\La}$ and the set of virtual arrows. We note that in the setting of this paper, the permutation $g$ will not have
3-cycles and the second exception in (ii) and in (iii) cannot occur. We also assume that $c_\theta=1$ for all 
virtual arrows $\theta$. \smallskip 

Moreover, note that this is the case for the spherical (and almost spherical) algebras defined in the 
introduction, where virtual arrows appear. For instance, observe that the commutativity relations $(S1ab)-(S4ab)$ 
are obtained from the relations (i) by substituting $\beta\nu,\delta\alpha,\omega\gamma$ and $\sigma\rho$ 
instead of the virtual arrows $\xi,\eta,\ve$ and $\mu$, respectively (and removing the relations 
$\theta f(\theta)- \bar{\theta}$, for all virtual arrows $\theta$). The zero relations $(Z1)-(Z16)$ are 
obtained from the relations (ii) and (iii) by analogous substitutions.

\subsection{Paths of length $\leq 3$}\label{2.2}

First we recall  the following result from \cite{EHS1}, which
is a main tool also to identify minimal relations. Throughout $\La$ is a TSP4 algebra.

\begin{lemma}\label{lem:2.1} Assume $\sigma\tau$ is a path $i\to j$ in the path algebra of $KQ^S$ where  $\sigma,  \tau$ are arrows. If $\sigma\tau$ is involved in some minimal relation then there must be an arrow from $j$ to $i$.
\end{lemma}

This implies that any of the  path of length two between 1-vertices  of $Q^{S}$ 
cannot  occur in a minimal relation of the algebra. However, the algebra
must be tame and this gives the following.

\begin{lemma}\label{lem:2.2}
	(a) \ Given two paths of length three in $Q^S$, both from $i$ to $j$.
Then there must be some minimal relation involving these paths.\\
(b) \ The same holds for the case of $Q^{S'}$ if the paths do not involve 
$\ve$ or $\mu$.
\end{lemma}

\begin{proof}

 We give the argument for such paths starting at a 1-vertex, 
say at $b_1$. Consider
the subquiver formed by $\beta\gamma\sigma$ and $\beta\nu\delta$. If there is no minimal relation involving the two paths of length three then this subquiver 
 is wild hereditary, and it follows that $\La$ is wild.
 A similar argument is valid for such paths ending at a 1-vertex. \end{proof}

\medskip

\begin{lemma}\label{lem:2.3}
Assume $\La$ has Gabriel quiver $Q^{S'}$.  Then we have:\\
(a) \ A path of length two which is not part of one of the triangles $(\sigma, \rho, \ve)$ or $(\mu, \omega, \gamma)$ 
does not occur in a minimal relation.\\
(b)\  Any path of length two which is part of one of the triangles must occur in some minimal relation. 
Then it belongs to $J^3$.
\end{lemma}

\begin{proof} Part (a) follows easily from Lemma \ref{lem:2.1}. For (b), we consider the idempotent algebra $e\La e$, 
where $e = 1- e_{d_2} - e_{b_1}$. This has a quiver $\wt{Q}$, whose vertex set is obtained from $Q_0$ by removing 
$b_1$ and $d_2$. Since $\alpha\beta,\nu\delta\nprec I$, the arrows of $Q$ around the square 
$(\alpha \ \beta \ \nu \ \delta)$ are changed, we only have $\wt{\nu} = \nu\delta:2\to 1$ and 
$\wt{\alpha} = \alpha\beta:1\to 2$ and arrows not around this block remain arrows for the idempotent algebra. 
Thus $\wt{Q}$ is $2$-regular. \smallskip 

Now, observe that $\rho\omega \nprec I$, so it follows from \cite[Proposition 4.2]{GQT} that $\rho\ve\prec I$, since 
$e\La e$ is a tame algebra with $2$-regular Gabriel quiver. Similarly, since $\ve\mu \nprec I$, we get $\ve\sigma\prec I$. 
Then $\sigma\rho \prec I$, by \cite[Lemma 4.3]{EHS1}. Arguments for the paths in the second triangle are analogous. 
\end{proof}

\subsection{Period four}\label{2.3}

Take a simple module $S_i$, we assume that it has $\Omega$-period four and 
$\Omega^3(S_i)\cong \Omega^{-1}S_i \cong P_i/S_i$ where $P_i = e_i\La$, the 
indecomposable projecte with top and socle isomorphic to $S_i$. This means that there
is an exact sequence, with projective modules $P$ and $P'$, of the form
$$0\to \Omega^{-1}(S_i)=P_i/S_i \to P' \to P \stackrel{d}\to {\rm rad}P_i = \Omega(S_i) \to 0$$
and  $\Omega^2(S_i)$ is isomorphic to the kernel of $d$. 
If $i$ is a 1-vertex and $\tau$ is the arrow starting at $i$
then  ${\rm rad} P_i = \tau \La$,  and $P_i/S_i = \tau^*\La$ where
$\tau^*$ is the arrow ending at $i$. 
If $i$ is a 2-vertex, ${\rm rad} P_i = \tau\La + \bar{\tau}\La$, 
and $P_i/S_i\cong (\tau_1^*, \tau_2^*)\La$ wher $\tau, \bar{\tau}$ start at $i$ and $\tau_1^*, \tau_2^*$ end at $i$. 
Then  the minimal relations starting at $i$ give rise to minimal generators of $\Omega^2(S_i)$. There is either  one such relation when $i$ is a 1-vertex, or
there are two.  
We spell this out for a 1-vertex in our setting, the case of a 2-vertex will
be written down later.

The radical of $e_{b_1}\La$ is isomorphic to $\beta\La$, and $e_{b_1}\La/S_{b_1}$ is isomorphic to $\alpha\La$. Hence the 
above exact sequence is of the form
	$$0 \to \alpha\La \to P_1 \to P_2 \to \beta\La \to 0$$
We have a minimal relation involving paths of length three from vertex $b_1$ 
	to vertex $1$, and this gives the following minimal relations:
$$r_1(\beta\nu\delta) + r_2(\beta\gamma\sigma) + \beta E=0 \ \ \mbox{and }\ \   r_1(\nu\delta\alpha) + r_2(\gamma\sigma\alpha) + E\alpha = 0$$
where $r_i\in K$ not both zero and $E\in e_2J^4e_1$ (in case $Q^S$, the error term $E$ is in $J^6$).

\section{\bf Minimal relations and indecomposable projectives}\label{3}

In this section, we desribe the general shape of relations and generators of projective modules. 

\subsection{Minimal relations}\label{3.1} 

As we have seen above, minimal relations involving paths of length three
starting or ending at $1$-vertices come in pairs, and using this for each 
1-vertex we obtain the following. 

\begin{proposition} \label{prop:3.1} Let $\La$ be a TSP4 algebra with 
	Gabriel quiver $Q^S$. 
Then we may assume that 
the  minimal relations involving paths of length three, contained in $I$, are as follows, up to labelling:
\begin{align*} 
(1a) &\  \beta\nu\delta -c (\beta\gamma\sigma) + \beta E= 0,&&&
(1b) &\  \nu\delta\alpha -c(\gamma\sigma\alpha)  + E\alpha = 0, &\\
(2a) &\  \delta\alpha\beta -d(\delta\rho\omega) + \delta F=0, \ &&& 
(2b) &\  \alpha\beta\nu -d(\rho\omega\nu) + F\nu=0, &\\
(3a) & \ \sigma\rho\omega -c(\sigma\alpha\beta) + \sigma G=0, &&&
(3b) & \ \rho\omega\gamma -c(\alpha\beta\gamma) + G\gamma=0, \\
(4a) & \ \omega\gamma\sigma -d(\omega\nu\delta) + \omega H=0, &&&
(4b) & \ \gamma\sigma\rho -d(\nu\delta\rho) + H\rho=0.
\end{align*}
Here $c, d\in K$ and $E, F, G, H$ are elements in $J^6$.\\

When $\La$ is a TSP4 algebra with Gabriel quiver $Q^{S'}$, we may assume as well that
 the  minimal relations (1a), (1b) and (2a), (2b) hold.
\end{proposition}

\begin{proof} Recall that we make the choice for the blocks, such that the paths of length three around a block 
should all occur in some minimal relations (Lemma \ref{lem:2.2}(a); see also \cite[..Square Lemma]{EHS1}). We 
exploit the exact sequences for the simples at 1-vertices. As explained above, generators of the second syzygies 
correspond to minimal relations. Take the simple module $S_{b_1}$, the identities are written down in
\ref{2.3}. We require $\beta\nu\delta$ to occur in a minimal relation
so $r_1\neq 0$ and we may take $r_1=1$ and $r_2=-c$.
Similarly exploiting the sequence
for $S_{d_2}$ gives identities (2a) and (2b), denoting the scalar by $-d$.
We may similarly write down identities (3ab) and (4ab) using the other two simple modules, with
scalars $-c_1$ and $-d_1$. 

In the case when the Gabriel quiver is $Q^S$,  we must show that $c=c_1$ and $d=d_1$. For this we use the 
exact sequence for $S_1$, it is of the form
$$0\to (\delta, \sigma)\La \to P_{d_2}\oplus P_{b_2} \to P_{b_1}\oplus P_{d_1} \to \alpha\La + \rho\La \to 0.
$$
We observe that the equations (2b) and (3b) give rise to two independent  minimal generators for $\Omega^2(S_1)$. At the moment
we do not need error terms explicitly, they are in $J^7$. Writing 
$F\nu = \alpha A + \rho B$, and  $G\gamma = \alpha A' + \rho B'$,  we get  generators
$$(\beta\nu + A, -d(\omega\nu) + B) \ \mbox{and} \ (-c_1(\beta\gamma)+A',\omega\gamma+B').$$
Let $M$ be the matrix whose columns are the   generators of $\Omega^2(S_1)$. Then we have $M{\delta\choose \sigma} = 0$, 
which gives the following two minimal relations 
$$\beta\nu\delta -c_1(\beta\gamma\sigma)  + A\delta+ A'\sigma=0, \mbox{ and }
\omega\gamma\sigma  -d(\omega\nu\delta) +  B'\sigma + B\delta=0.$$ 
Now, using equations (1a) and (4a), we deduce that $c_1=c$ and $d_1=d$. This finishes the proof. 
\end{proof} 

\medskip

This Propsition gives explicit minimal relations for paths 
of length three around a block of type V. We will therefore use the paths which are not around a block of type V (and not around a triangle, in the case of $Q^{S'}$)
towards finding bases for the projective modules. 
We set
$$X_{\alpha}:= \alpha\beta\gamma\sigma, \  \ \mbox{ and} \ \ Y_{\rho}:= \rho\omega\nu\delta.
$$
We use similar notation for the rotations of these paths, for example $X_{\beta} = \beta\gamma\sigma\alpha$.

\subsection{Spanning sets for projectives}\label{3.2}

(1) \ By using the minimal relations for the paths of length three, and the information on paths of length two,
we may  determine inductively spanning sets for the radical quotients of 
$e_1\La$ and $e_2\La$. For the case of $Q^{S'}$, we use that paths of length two around a triangle lie in $J^3$, see Lemma \ref{lem:2.3}(b).  
We consider the case of  $e_1\La$ in some detail, the case of $e_2\La$ is analogous.
The first radical quotient  has basis formed by the cosets of $e_1, \alpha$ and $\rho$.
Next, 
$e_1J^2 = \langle \alpha\beta, \rho\omega \rangle + e_1J^3$
and 
$e_1J^3$
is generated by $\alpha\beta\gamma, \ \alpha\beta\nu, \ \rho\omega\nu, \ \rho\omega\gamma$. 

Consider $\alpha\beta\nu$, this is either in $J^7$ (if $d=0$), or else 
$\alpha\beta\nu\equiv \rho\omega\nu$. We keep $\rho\omega\nu$ as a generator
and do not need $\alpha\beta\nu$. 
Consider $\rho\omega\gamma$. This is either in $J^7$ (if $c=0$), or else 
$\rho\omega\gamma\equiv \alpha\beta\gamma$. We keep $\alpha\beta\gamma$ as a generator, 
and we do not need $\rho\omega\gamma$. 
This shows that  $e_1J^3$ is generated by $\alpha\beta\gamma$ and $\rho\omega\nu$, and then 
$e_1J^4$ is generated $X_{\alpha}$ and $Y_{\rho}$ modulo $J^5$. 
Iterating this, we obtain a spanning set, which consists of all initial submonomials of $X_\alpha^n$ and 
$Y_\rho^n$, $n\in\bN$.
We can now identify socle elements for $e_1\La$. Let $m, m'$ be maximal such
that $W_{\alpha}:= X_{\alpha}^m\neq 0$, and $W_{\rho} := Y_{\rho}^{m'}\neq 0$.
Using that $\La$ is symmetric, they must span the socle of $e_i\La$,  and  $W_{\alpha} \equiv W_{\rho}$.
This shows that $e_1\La$ is spanned by  the initial
submonomials of $W_{\alpha}$ and  $W_{\rho}$.
We write the initial submonomial of length $t$ of $W_{\alpha}$ as
$[W_{\alpha}]_t$, then the spanning set is  
$$\{ [W_{\alpha}]_t \mid 0\leq t\leq 4m\} \cup \{ [W_{\rho}]_s \mid 1\leq s < 4m'\}.
$$
There is an analogous spanning set for $e_2\La$.

(2) \ Taking rotations of $W_{\alpha}$ and $W_{\rho}$ gives  socle elements for the other
indecomposable projective modules. We denote the socle element starting with $\tau$ by $W_{\tau}$, for $\tau$ any arrow.
Furthermore, we denote by $A_{\tau}$ the initial submonomial of $W_{\tau}$ such that $A_{\tau}\tau^*= W_{\tau}$ for the arrow $\tau^*$ with $g(\tau^*)=\tau$. 
We will see later that $A_{\tau}$  is the element occuring in the presentation
for the WSA.

\bigskip
\begin{corollary}\label{cor:3.2} Let $E, F, G$ and $H$ be elements in the identities of Proposition \ref{prop:3.1}. 
We can take $E$ and $G$ to be linear combinations of powers of $X_{\alpha}$ and $X_{\gamma}$; and we can take $F$ and 
$H$ to be linear combinations of powers of $Y_{\rho}$ and $Y_{\nu}$. \end{corollary}

\begin{proof}

 Consider for example the element $E$, it occurs as
$\beta E$ and $E\alpha$, so we can take it  as stated, by using the above spanning set. \end{proof}

\medskip

\begin{lemma}\label{lem:3.3} 
(a) Assume $c\neq 0$ and $d=0$, then $m=1$. Similarly if $c=0$ and $d\neq 0$, then $m'=1$.\\
(b) If $d=0$ then $m'> 1$ and if $c=0$ then $m >1$.\\
(c) If $c, d$ are both non-zero then $X_{\alpha}\equiv Y_{\rho}$. 
For the quiver $Q^{S'}$, this implies that  $m=1=m'$.
\end{lemma}

\begin{proof}

  (a) 
 Let $c\neq 0$ and $d=0$, then by (1a) and (2b) of Proposition \ref{prop:3.1},  we have
 $$(\alpha\beta\gamma\sigma)\alpha  \equiv (\alpha\beta\nu\delta)\alpha  \equiv  -F\nu\delta\alpha \in J^9
$$
Hence it lies in $J^9\cap e_1\La e_{b_1}$ which is the span of the set 
of all $X_{\alpha}^j\alpha$. But then 
$X_{\alpha}\alpha = \sum c_jX_{\alpha}^j\alpha$ which implies that the product of $X_{\alpha}\alpha$ and a unit 
is zero, and hence $X_{\alpha}\alpha=0$.
Therefore $X_{\alpha}=W_{\alpha}$ and $m=1$. The other part of (a) is similar.

(b) Assume $d=0$, and assume for a contradiction that $m'= 1$. 
Then $F=0$ and  identity (2a) gives $\delta\alpha\beta =0$. Then $\delta\alpha J=0$ and $\delta\alpha$ is in the socle. 
But it belongs to $e_{d_2}\La e_{b_1}$ and is therefore zero, a contradiction to Lemma \ref{lem:2.1}.
Similarly one gets the other case.

(c) 
Assume $c, d$ are both non-zero,  then  by (1a) and (2b) we have $X_{\alpha} \equiv Y_{\rho}$. Assume we have
$Q^{S'}$. Then $\sigma\rho\in J^3$ (see Lemma \ref{lem:2.3}(b)) and it follows that
$Y_{\rho}\rho\in J^6$. It also lies in $e_1\La e_{d_1}$ and then it is a linear combination of elements $Y_{\rho}^j\rho$. 
The argument used in part (a) of this proof shows that $Y_{\rho}\rho =0$ and $m'=1$. Since $X_{\alpha}\equiv Y_{\rho}$ 
we also have $m=1$. \end{proof}

\bigskip

Note that in case $m$ or $m'$ is $1$, the proper initial submonomials of $W_\alpha$ and $W_\rho$ end at different 
vertices (except $2$), so these are linearly independent, since we cannot have a commutativity relations involving 
paths $\alpha\beta,\rho\omega$. Hence Lemma \ref{lem:3.3} shows that the above spanning 
sets for $e_1\La$ and $e_2\La$ are bases, unless possibly $c, d$ are both non-zero and the algebra has Gabriel quiver $Q^S$. 
We consider now separate cases.

\section{\bf The proof in case $c$ or $d$ is zero, or the quiver is $Q^{S'}$ }\label{4}

In this case, the spanning sets in \ref{3.2} are bases, and we have the conditions on $m, m'$ as in Lemma 
\ref{lem:3.3}. With these, we can now determine bases for the projective modules at 1-vertices.

\medskip

\begin{lemma} \label{lem:4.1} Assume at least one of $c$ or $d$ is zero, then the projectives 
$P_{b_1}$, $P_{b_2}$, and $P_{d_1}$, $P_{d_2}$ have bases 
	$$\{ [W_{\beta}]_t \mid 0\leq t\leq 4m \} \cup \{ \beta\nu\}, \ \  \{ [W_{\sigma}]_t \mid 0\leq t\leq 4m\}\cup \{ \sigma\rho\}$$
and
$$
\{ [W_{\delta}]_t\mid 0\leq t\leq 4m'\} \cup \{ \delta\alpha\}, \ \ \{ [W_{\omega}]_t \mid 0\leq t\leq 4m'\}\cup \{ \omega\gamma\},$$ 
respectively. For the $Q^{S'}$ case, the projectives  $P_{b_1}$ and $P_{d_2}$ have the same bases.\end{lemma}

\begin{proof} The proof that this is a spanning set follows from similar arguments as in \ref{3.2}. This 
is elementary exercise to show that initial submonomials of a given cycle in the socle are linearly independent. 
For example, paths in $\{[W_\beta]_t\mid 0 \leq t 4m\}$ are independent. Since $\beta\nu$ ends at different 
vertices that the submonomials of $W_\beta$, we get that the indicated set is a basis. Similarily in other 
cases. \end{proof} 

\medskip

Recall our convention about socle elements and maximal subpaths, so that
for example $W_{\beta} = A_{\beta}\alpha$.

\medskip

\begin{corollary}\label{cor:4.2} We have the following zero zerlations 
$A_{\beta}\rho = 0$, $A_{\delta}\gamma=0$, $A_{\omega}\alpha=0$ and $A_{\sigma}\nu=0$. \end{corollary}

\begin{proof} By considering  the basis of $e_{b_1}\La$ we see that $A_{\beta}$ spans $e_{b_1}J^{4m-1}$. 
Therefore  $A_{\beta}\rho \in e_{b_1}J^{4m}$ which is spanned by $W_{\beta}$. Hence  
$A_{\beta}\rho = \lambda W_{\beta}$ for a scalar $\lambda$. Postmultiplying  with $e_{b_1}$ shows that $\lambda=0$. 
Similarly the other identities hold.
\end{proof}

We have another consequence in case $Q^S$. 

\begin{lemma}\label{lem:4.3} (a) \ The spaces 
	$e_{b_1}\La e_{d_1}$ and $e_{d_2}\La e_{b_2}$ are zero. Moreover
	$e_{b_1}\La e_{d_2}$ and $e_{d_1}\La e_{b_2}$ are 1-dimensional.\\
	(b) All paths of length three around a block of type V$_2$ belong
	to the second socle of the algebra. 
\end{lemma}

\begin{proof}
 (a)  The dimensions follow directly.\\
(b) \ 
We give details in one of the cases. Consider the path $\beta\nu\delta$. 
Using the basis for $P_{b_1}$ from Lemma \ref{lem:4.1}, we conclude that $\beta\nu\delta$ is a linear combination 
of paths $[W_\beta]_t$, $t\geqslant 0$, which end at $1$. Hence, we have the following relation in $\La$ 
$$\beta\nu\delta=[W_\beta]_t(1+C),$$ 
where $[W_\beta]_t$ is an initial submonomial of $W_\beta$ of the smallest length in this combination, 
and $C\in e_1 J e_1$. As a result, $1+C$ is a unit, and therefore, we may write  
$$\beta\nu\delta=c_\beta (\beta\gamma\sigma\alpha)^{m_1-1}\beta\gamma\sigma,$$ 
for some $m_1\in\{1,\dots,m\}$ and a scalar $c_\beta\in K^*$. Now, it follows that $m_1=m$, since otherwise 
we get $\beta\nu\delta\alpha\beta\gamma=c_\beta X_\beta^{m_1}\beta\gamma\neq 0$, so 
$\delta\alpha\beta\gamma$ is a non-zero element in $e_{d_2}\La e_{b_2}$, a contradiction with (a). Therefore, 
indeed $\beta\nu\delta=c_\beta A_\beta$ belongs to the second socle (use also $\beta\nu\delta\rho\in e_{b_1}\La e_{d_1}=0$) 
and $\beta\nu\delta\alpha=c_\beta W_\beta$ is in the socle. Involving an analogous basis for the  projective left module 
$\La e_{d_2}$ one can show that $\alpha\beta\nu=c_\rho A_\rho$ is in the second socle and $\alpha\beta\nu\delta$ 
is in the socle. The proof for the remaining two rotations is similar. The same arguments work for the paths 
in the second block V$_2$ (in case $Q^{S'}$). 
\end{proof}

We only mention that (b) holds also in case $Q^{S'}$, but it needs the description of bases given later 
(see Lemma \ref{lem:6.1}).  

\bigskip

The following proves part of Theorem \ref{thm:1.1}.

\begin{proposition} \label{prop:4.4} 
Assume that the Gabriel quiver of $\La$ is of the form 
	$Q^S$, and that at least one of $c, d$ is zero. Then
the algebra $\La$ is isomorphic to a WSA, with virtual arrows $\xi, \eta$ and $\mu, \ve$, and $f$. 
The multiplicities are $m_{\alpha} =m$ and $m_{\rho} = m'$, and $m_{\xi}=1=m_{\mu}$. 
\end{proposition}

\begin{proof} Without loss of generality, we may assume that $c=0$, so $m>1$, by Lemma \ref{lem:3.3}(b). Applying 
Lemma \ref{lem:4.3}(b), we conclude that all paths of length $3$ around a block of 
type V$_2$ are in the second socle. In particular, using bases of projectives, the  
relations $(1a)-(4a)$ become 
$$\beta\nu\delta=c_\beta(\beta\gamma\sigma\alpha)^{m-1}\beta\gamma\sigma, \ 
\delta\alpha\beta=c_\delta(\delta\rho\omega\nu)^{m'-1}\delta\rho\omega, \ 
\sigma\rho\omega=c_\sigma(\sigma\alpha\beta\gamma)^{m-1}\sigma\alpha\beta, \ 
\omega\gamma\sigma=c_\omega(\omega\nu\delta\rho)^{m'-1}\omega\nu\delta,$$ 
for some scalars $c_\beta,c_\delta,c_\sigma,c_\omega\in K^*$. In other words, we have 
$$\theta f(\theta) g(f(\theta))=c_\theta A_\theta,$$ 
for any arrow $\theta$ starting from a $1$-vertex (note that $f$ is always defined, and it is induced from $Q$). 
As a result, we conclude from Proposition \ref{prop:3.1} that the relations $(1b)-(4b)$ admit analogous form 
$$\theta g(\theta) f(g(\theta))=c_{\bar{\theta}}A_{\bar{\theta}},$$ 
where $c_\gamma=c_\beta$, $c_\rho=c_\delta$, $c_\alpha=c_\sigma$ and $c_\nu=c_\omega$. Now, taking the further 
relations $(2b)$ and $(3b)$ gives a matrix $M$, whose columns are generators of $\Omega^2(S_1)$. 
As in the proof of Proposition \ref{prop:3.1}, we use the exact sequence for $S_1$, and obtain two other 
relations from the equality $M\cdot\vec{\delta \\ \sigma}=0$. These compared with $(1a)$ and $(4a)$ give  
$c_\alpha=c_\beta$ and $c_\rho=c_\omega$. \smallskip 

This shows that $c_\bullet$ gives rise to a function constant on $g$-orbits, and consequently, the relations 
$(1ab)-(4ab)$ are exactly the relations $(S1ab)-(S4ab)$ defining the spherical algebra (with scalars $a=c_\beta$ 
and $b=c_\rho$). Moreover, it is easy to see that relations $(Z1)-(Z16)$ defining the spherical algebra 
are also satisfied in $\La$, because the paths of length four around blocks V$_2$ are in the socle, by 
Lemma \ref{lem:4.3}(b), and the remaining relations follow from the vanishing conditions in Lemma \ref{lem:4.3}(a). \medskip 

Concluding, if $\hat{\La}$ is the spherical algebra $\hat{\La}=KQ^S/I^S$ given by weights $m_\alpha=m$, 
$m_\rho=m'$, and scalars $a=c_\beta$, $b=c_\rho$, then $\La$ satisfies all relations defining $\hat{\La}$. 
Therefore, we have an epimorphism $\hat{\La}\to\La$. Comparing known bases (hence dimensions) of the 
indecomposable projective $\La$-modules with those for $\hat{\La}$, we conclude that $\psi$ is an isomorphism, 
and the proof is now finished. 
\end{proof}

\bigskip

\section{Algebras when all scalars are non-zero}\label{5}

We will show that the algebra is either a weighted surface algebra
with all  multiplicities equal to 1, in the same family as the
algebras in the previous case,  or else it is a Higher Spherical Algebra.
That is, it isomorphic to an algebra denoted by $S(m, \lambda)$ in \cite{HSA},
with Gabriel quiver $Q^S$ subject to the following relations, with
 elements $A_{\alpha}$ as defined before, and $\lambda\in K^*$:
\begin{align*}
	 &&\  \beta\nu\delta = \beta\gamma\sigma + \lambda A_{\beta},   &&  &&\  \alpha\beta\nu  =  \rho\omega\nu,  &&\\
	 &&\  \nu\delta\alpha = \gamma\sigma\alpha +\lambda A_{\gamma}, &&  &&\  \delta\alpha\beta = \delta\rho\omega,  &&\\
	&& \ \sigma\rho\omega  = \sigma\alpha\beta + \lambda A_{\sigma}, &&  && \ \omega\gamma\sigma  = \omega\nu\delta,  &&\\
	 && \ \rho\omega\gamma  = \alpha\beta\gamma + \lambda A_{\alpha}, &&  && \ \gamma\sigma\rho = \nu\delta\rho.&&\\
	&& \ 	(\alpha\beta\gamma\sigma)^m\alpha=0,  &&&& (\gamma\sigma\alpha\beta)^m\gamma = 0. 
\end{align*}

\subsection{Bases for the projective modules}\label{5.1}

We had already found spanning sets for $e_1\La$ and $e_2\La$, and we determine now a basis.
Let as before
$$X:= X_{\alpha} = \alpha\beta\gamma\sigma \ \  \mbox{and} \ \ Y:= Y_{\rho} = \rho\omega\nu\delta.
$$
In the following we will use the identities of Proposition \ref{prop:3.1}.
By (3b) and (4a) we have
$$X \equiv \rho\omega\gamma\sigma \equiv \rho\omega\nu\delta = Y,$$
hence $X^r\equiv Y^r$ for $r\geq 1$.
Moreover $\alpha\beta\nu\delta \equiv \rho\omega\nu\delta$ (by (2b)) and
$\rho\omega\gamma\sigma \equiv \alpha\beta\gamma\sigma$ (by (3b)).
That is, all paths of length four starting at vertex $1$  are linearly
dependent modulo $J^5$.
In particular, the local algebra at vertex $1$ is generated by $X$ and has finite type.

We will use rotations of $X$ and of $Y$, to construct bases for the indecomposable projective modules.
We write $X_{\beta} = \beta\gamma\sigma\alpha$ for the rotation of $X$,
similarly we use
the notation $X_{\gamma}$ and $X_{\sigma}$. Similarly we write rotations
of $Y$.
With these, the following describes bases for $e_1\La$ and also for $e_2\La$, the proof is straightforward.

\bigskip

\begin{lemma} \label{lem:5.1}
(a) The module $e_1\La$ has basis
$$\{ [X^m]_i \mid 0\leq i\leq 4m \} \cup \{   X^i\rho, X^i\rho\omega\nu\mid 0\leq i\leq m-1\}\cup\{\rho\omega\}$$
There is also the $Y$-basis, obtained from the $X$-basis by replacing $X$ by $Y$, and $\rho, \omega, \nu$ by $\alpha, \beta, \gamma$.

(b) The module $e_2\La$ has basis
$$\{ [X_{\gamma}^m]_i \mid 0\leq i\leq 4m \} \cup \{   X_{\gamma}^i\nu,  X_{\gamma}^i\nu\delta\rho \mid 0\leq i\leq m-1\} 
\cup\{\nu\delta\}$$
There is also the $Y$-basis, obtained from the $X$-basis by
replacing $X$ by $Y$ and $\nu, \delta, \rho$ by $\gamma, \sigma, \alpha$.
\end{lemma}

\bigskip
The dimension vector of $e_1\La$  is $(m+1, m, m+1, m, m, m)$ using the order $1, b_1, 2, b_2, d_1, d_2$, and the dimension vector
for $e_2\La$ is the same.
We identify a basis for $P_{b_1}$. The socle of $P_{b_1}$ is  spanned by $X_{\beta}^m$, and  we have a basis
$$\{ [X_{\beta}^m]_t, 0\leq t\leq 4m\} \cup \{  X_{\beta}^r\beta\nu, 
\ X_{\beta}^r\beta\nu\delta\rho,   0\leq j\leq m-2\} \cup \{ X_{\beta}^{m-1}\beta\nu\}.$$
There are similar bases  for the  projectives at the other 1-vertices.

 \subsection{Changing generators}\label{5.2}

As before, we use the identities from Proposition \ref{prop:3.1}. 
We will change generators,  so that identities (2ab) and (4ab) become the four relations with only two terms   in the presentation for the HSA.
We start by refining the error terms $F$ and $H$.
The element $F\in e_1\La e_2$ is multiplied with $\delta$ and $\nu$, this means that we
should express it in terms of the $Y$-basis.
We write
it as
$$F=  \rho F_{d_1}\omega$$
where $F_{d_1}$ is a linear combination of powers of $Y_{\omega}$.
Similarly we write
$$H=\nu H_{d_2}\delta.
$$
With these, identities (2a) and (2b) become 
$$\delta\alpha\beta = \delta\rho(d - F_{d_1})\omega, \ \ \alpha\beta\nu = \rho(d -  F_{d_1})\omega\nu.
$$
We also rewrite the identities (4a) and (4b), this gives
$$\omega\gamma\sigma = \omega\nu(d - H_{d_2})\delta, \ \ \  \gamma\sigma\rho = \nu(d - H_{d_2})\delta\rho
$$
We replace $\rho$ by $\rho'$ where
$\rho':= \rho(d-F_{d_1})$.
Then (2a) and (2b) become
$$\delta\alpha\beta = \delta\rho'\omega \ \mbox{and} \ \alpha\beta\nu = \rho'\omega\nu \leqno{(5.2.1)}$$
We also replace $\nu$ by  $\nu'$ where
$\nu' = \nu(d-H_{d_2}).$
 Then (4a) and (4b) become
 $$\omega\gamma\sigma = \omega\nu'\delta \ \ \ \mbox{and} \ \ \gamma\sigma\rho = \nu'\delta\rho. \leqno{(5.2.2)}
 $$
We postmultiply the second identity in (5.2.1) by the unit $(d - H_{d_2})$ so that we get 
$$\alpha\beta\nu' = \rho'\omega\nu'$$
Similarly postmultiplying the second identity in (5.2.2) by the unit $(d-F_{d_1})$, we get
$$\gamma\sigma\rho' =\nu'\delta\rho'.$$

To complete the work, 
we must rewrite the identities (1ab) and (3ab). We substitute
$$\rho= d^{-1}\rho' + d^{-1}\rho F_{d_1}, \ \ \nu = d^{-1}\nu' + d^{-1}\nu H_{d_2}$$
and incorporate the terms with $F$ and $H$ into $E$ and $G$, 
then we multiply the identities by $d^{-1}$, which gives a scalar factor $c_1:= d^{-1}c$. We keep the terms of length three, with scalar $c_1$, and we
write the identities again, with new error terms $E'$ and $G'$, note that the identies come in pairs.
\begin{align*} (1ab) & \beta \nu'\delta - c_1(\beta\gamma\sigma) + \beta E' = 0 \ \ \mbox{and}  &  
\nu'\delta\alpha - c_1(\gamma\sigma\alpha + E'\alpha=0 \cr
	(3ab) &  \sigma\rho'\omega - c_1(\sigma\alpha\beta) + \sigma G'=0 \ \ \mbox{and}&  \rho'\omega\gamma - c_1(\alpha\beta\gamma) + G'\gamma = 0
\end{align*}
We may write $E'$ and $G'$ in the $X$-basis as $E'=\gamma E''\sigma$, and  $G' = \alpha G''\beta$.

\begin{lemma} We have $\beta\gamma E''\sigma = G''\beta\gamma\sigma$.
\end{lemma}

\begin{proof}

 We use the exact sequence for $S_1$. The generators of $\Omega^2(S_1)$ can be taken from using identities (2a) and (3b). Writing them as
a matrix this is
$$M= \left(\begin{matrix} \beta\nu' & -c_1\beta\gamma + G''\beta\gamma\cr -\omega\nu' & \omega\gamma\end{matrix}\right).
$$
Then $M{\delta \choose \sigma} = 0$ gives two identities, the second is (4a) and the first one has the same length three terms as (1a).
The length three terms are the same, and the claim follows.
\end{proof}

Now we have four identities,  all error terms are expressed in the X-basis.
Next we obtain a zero relation:

 \begin{lemma} \label{lem:5.3} Let $Y=\rho'\omega\gamma\sigma$. Then
 $$Y^2  = c_1Y^2 - \rho'\omega\gamma\sigma G'\gamma\sigma.$$
 Hence if $c_1\neq 1$ then $Y^2 =0$.
 In any case we have
 $\rho'\omega\gamma\sigma G'\gamma\sigma =0$.
 \end{lemma}

 \begin{proof}
 We have
 $$\begin{aligned} Y^2 =  \rho'\omega\gamma\sigma\rho'\omega\gamma\sigma 
  =_{(3a)}    \rho'\omega\gamma(c_1\sigma \alpha\beta \ - \ \sigma G')\gamma\sigma \cr
  =_{(4a)}  c_1\rho'\omega\nu'\delta \alpha\beta\gamma\sigma -  \rho'\omega\gamma\sigma G'\gamma\sigma
    =_{(2a)}  c_1\rho'\omega\nu'\delta \rho'\omega\gamma\sigma -  \rho'\omega\gamma\sigma G'\gamma\sigma\cr
     =_{(4a)}  c_1\rho'\omega\gamma\sigma \rho'\omega\gamma\sigma -  \rho'\omega\gamma\sigma G'\gamma\sigma 
	 =  c_1Y^2 - \rho'\omega\gamma\sigma G'\gamma\sigma. &&
 \end{aligned}
  $$
  The claim follows.
   \end{proof}
    \bigskip

 \begin{corollary} The element $G'\gamma\sigma$ is in the socle, and  
$G'\gamma$ and $\sigma G'$ belong to the second socle.
The element $\alpha\beta E'$ is in the socle, and $\beta E'$ and $E'\alpha$ belong to the second socle.
 They are non-zero.
 \end{corollary}

\begin{proof} By identity (3b) we have
$$G'\gamma\sigma = c_1\alpha\beta\gamma\sigma - \rho'\omega\gamma\sigma$$
note that the last term is the element $Y$ from Lemma \ref{lem:5.3}.
The element $G'\gamma\sigma$ is in $e_1\La e_1$ and it is annihilated by the generator of the local algebra
 which can be taken as $\rho'\omega\gamma\sigma$, by Lemma \ref{lem:5.3}. Hence $G'\gamma\sigma$ is in the socle of the local algebra.
 Now, $G'\gamma$ is in $e_1\La e_{b_2}$ and $\sigma$ is the only arrow starting at vertex $b_2$. Hence $G'\gamma$ is in the second socle.
By rotation, also $\sigma G'\gamma$ is in the socle, and the previous argument
 shows that $\sigma G$ is in the second socle.

Next, $\beta E' = \beta\gamma E''\sigma = G''\beta\gamma\sigma$ and hence
$\alpha\beta E' = G'\gamma\sigma$ which is in the socle.
Therefore $\beta E'$ is in the second socle. The rotation $\beta E'\alpha$ is in the socle and $E'\alpha$ is in the second socle.
$\Box$

\bigskip

When $Y^2=0$ we have $m=1$, and the algebra is a WSA, by the argument
as in Proposition \ref{prop:4.4}.  Note however that for $m=1$, we must have $c_1\neq 1$ 
(otherwise the algebra is not symmetric, see \cite{WSA-GV}). 

We continue with $Y^2\neq 0$, so that $m>1$. We have
already obtained almost all minimal relations for a HSA, 
the nilpotence relation follows as well since $X_{\alpha}^m$
is non-zero in the socle, and so are all rotations.
We may compare dimensions and conclude that $\La$ is
a HSA. \end{proof}

\bigskip

\section{Almost spherical algebras}\label{6}

Assume $\La$ is  a TSP4 algebra with Gabriel quiver $Q^{S'}$.
As before, we will determine minimal relations, and then use this
to show it is isomorphic to a WSA.

Recall that identities (1ab) and (2ab) of Proposition \ref{prop:3.1} hold in this case, and we also
have determined bases for the projective modules except for those at $d_1$ and $b_2$.
Note that we already know that  $e_{b_1}\La e_{d_1}=0$ and that $e_{d_1}\La e_{b_1}=0$, since $\La$ is symmetric 
(see also Lemma \ref{lem:4.3}(a)).   

\bigskip

\begin{lemma} \label{lem:6.1} (a) \ The projective $e_{d_1}\La$ has basis $\{ [W_{\omega}]_t\mid 0\leq t \leq 4m' \} \cup \{ [W_{\ve}]_s\mid 1\leq s < 2m_{\ve}\}$ where $W_{\ve} = (\ve\mu)^{m_{\ve}}$
spans the socle of $e_{d_1}\La$.  \\
(b) The projective $e_{b_2}\La$ has basis $\{ [W_{\sigma}]_t\mid 0\leq t \leq 4m\}\cup \{ [W_{\mu}]_s\mid 1\leq s < 2m_{\ve}\}$ and $W_{\mu}$ is a rotation of $W_{\ve}$. \\
(c) \ The local algebra  at vertex  $d _1$ has independent generators $Y_{\omega}, Z_{\ve}:= \ve\mu$ and their products are zero. Similarly the local algebra at $b_2$ has independent generators  
$X_{\sigma}$ and $Z_{\mu}:= \mu\ve$, and their product is zero.
\end{lemma}

\begin{proof} Parts (a) and (b) are similar to \ref{3.1}. For part (c), for example $Y_{\omega}Z_{\ve}$ has a 
subpath $\delta\rho\ve \in e_{d_2}\La e_{b_2}=0$. 
\end{proof}

\bigskip

The following shows  that the remaining paths of length three
satisfy the essentially the  same identities as in the spherical case. 

\begin{lemma}\label{lem:6.2} Assume $\La$ is a TSP4 algebra with almost spherical quiver.
	Then we have
\begin{align*}
	(3a) & \ (\mu\omega) -c(\sigma\alpha\beta) + \sigma G=0, &&&
(3b) & \ (\rho\ve) -c(\alpha\beta\gamma) + G\gamma=0\\
(4a) & \ (\ve\sigma)  -d(\omega\nu\delta) + \omega H=0, &&&
(4b) & \ (\gamma\mu)  -d(\nu\delta\rho) + H\rho=0.
\end{align*}
Here $G, H$ are in $J^6$ and $c, d$ are as before.
\end{lemma}

\begin{proof} This is similar to previous proofs, using the 
the exact sequence for $S_{b_2}$,
$$0\to (\gamma, \ve)\La \to P_2\oplus P_{d_1} \to P_{d_1}\oplus P_1 \to \mu\La + \sigma\La \to 0$$
We omit the details.
\end{proof}

\bigskip
The following proves Theorem \ref{thm:1.2}.

\begin{proposition} Assume $\La$ is a  TSP4  algebra  with Gabriel quiver $Q^{S'}$. Then $\La$ is isomorphic to a WSA with data $(Q, f, m_{\bullet}, c_{\bullet})$ and where $\xi: b_1\to d_2$  and $\eta: d_2\to b_1$ are virtual arrows.
\end{proposition}

\begin{proof} This  uses a similar strategy as in Proposition \ref{prop:4.4}, but we show that the 
relations $(1ab)-(4ab)$ coincide with the relations $(S'1ab)-(S'4ab)$ defining a WSA given by the 
almost spherical quiver. The details for the relations involving paths around the triangles are 
left to the reader.  
\end{proof}


\begin{thebibliography}{99}
	

\bibitem{ASS}
  {I.~Assem, D.~Simson, A.~Skowro\'nski},
  {Elements of the Representation Theory of Associative Algebras 1:
  Techniques of Representation Theory},
  {London Mathematical Society Student Texts, vol. {65}},
  Cambridge University Press, Cambridge, 2006.


\bibitem{EHS1} K. Erdmann, A. Hajduk,  A. Skowyrski, {\it Tame algebras of period four}, Arch. Math. {\bf 122}(2024), 249-264




\bibitem{EHS2} K. Erdmann, A. Hajduk,  A. Skowyrski, {\it Local structure of tame symmetric algebras of period four} arxiv 2411.01235v1. 


\bibitem{EHS3} K. Erdmann, A. Hajduk,  A. Skowyrski, {\it Generalized quaternion type: biregular case}, 
to appear on arXiv. 


\bibitem{WSA}
K. Erdmann, A. Skowro\'{n}ski,  {\it  Weighted surface algebras.}
J. Algebra {\bf 505} (2018), 490-558.

\bibitem{GQT}
K. Erdmann, A. Skowro\'{n}ski,  {\it Algebras
of generalized quaternion type.}
                Adv. Math. {\bf 349}(2019), 1036-1116,

	\bibitem{WSA-GV} K. Erdmann, A. Skowro\'{n}ski,   {\it  Weighted surface algebras: general version.}
J. Algebra {\bf 544} (2020), 170-227.



\bibitem{WSA-GV-corr}
K. Erdmann, A. Skowro\'{n}ski,  {\it Weighted surface algebras: general version, Corrigendum.}
                J. Algebra {\bf 569} (2021), 875-889.

      
\bibitem{Sk}
A. Skowro\'{n}ski,
                {\it Selfinjective algebras: \ finite and tame type. } Contemp. Math. {\bf 406}(2006), 169-238.



\bibitem{HTA}
K.~Erdmann, A.~Skowro\'nski, 
{\it Higher tetrahedral algebras, }
Algebr. Represent. Theory {\bf 22}(2019), 387-406.
2018, 

\bibitem{HSA}
	K.~Erdmann, A.~Skowro\'nski, 
{\it Higher spherical algebras. } Arch. Math. {\bf 114} (2020), 25-39.


\bibitem{HSS} T. Holm, A Skowro\'{n}ski, A. Skowyrski, {\it Virtual mutations of weighted surface algebras.} J. Algebra {\bf 619} (2023), 822–859. 

\bibitem{SS} A. Skowro\'{n}ski, A. Skowyrski, {\it
	Generalized weighted surface algebras, } arXiv 2106.15218.

\bibitem{HT} A. Skowyrski, {\it Two tilts of higher spherical algebras}. Algebr. Represent. Theory {\bf 25}(2022), no. 1, 237-254.



\end{thebibliography}
\end{document}